\documentclass[a4paper,12pt]{article}
\usepackage{amssymb}
\usepackage{amsthm}
\theoremstyle{plain}
\newtheorem{Lemma}{Lemma}
\newtheorem{Proposition}{Proposition}
\newtheorem{Theorem}{Theorem}

\theoremstyle{definition}

\newcommand{\Sec}{\mathrm{Sec}}
\newcommand{\expdim}{\mathrm{expdim}}

\newcommand{\PP}{\mathbb{P}}
\newcommand{\C}{\mathbb{C}}

\begin{document}
\title{A Terracini Lemma for osculating spaces \\ 
with applications to Veronese surfaces}
\author{Edoardo Ballico and Claudio Fontanari}
\date{}
\maketitle 

\begin{small}
\begin{center}
\textbf{Abstract}
\end{center}
Here we present a partial generalization to higher order osculating spaces 
of the classical Lemma of Terracini on ordinary tangent spaces. 
As an application, we investigate the secant varieties to the osculating 
varieties to the Veronese embeddings of the projective plane. 
\vspace{0.5cm}


\end{small}

\vspace{0.5cm}

\section{Introduction}
The geometry of defective varieties, whose secant varieties have dimension 
less than the expected one, is a subtle and intriguing subject which has 
been investigated by several authors, both classical and modern (see for 
instance the introduction and the references in \cite{BF3}). The main 
tool for understanding defective varieties is provided by the celebrated 
Lemma of Terracini (see \cite{Terracini:11} for the original statement and 
\cite{Dale:84}, \cite{Adl:87} for modern versions): 

\begin{Lemma}\label{basicone} Let $X \subset \PP^r$ be an integral 
non-degenerate projective variety of dimension $n$ and let $h \le r$ be 
a positive integer. Take $h+1$ general points $p_0, \ldots, p_h$ of $X$ 
and let $P \in < p_0, \ldots, p_h >$ be a general point in the 
$h$-secant variety $\Sec_h(X)$ of $X$. Then the tangent space 
$T_P(\Sec_h(X))$ to $\Sec_h(X)$ at $P$ is given by $T_P(\Sec_h(X)) 
= < \bigcup_{i=0}^h T_{p_i}(X) >$. 
\end{Lemma}

The generalization of the above result from tangent spaces to higher 
order osculating spaces (respectively, from double points to points of 
higher multeplicity) is a very delicate problem, which has recently 
attracted a great deal of interest (for instance, we are aware of 
work in progress by Ciliberto on these topics). Here we present our 
attempts in this direction, which have been inspired by the less known 
paper \cite{Terracini:22} of Terracini: indeed, our Lemma~\ref{Terracini} 
should be regarded as a rough extension of Lemma~\ref{basicone}. 

As an application, we are going to investigate the secant varieties to the 
osculating varieties $T(m,V_{2,d})$ of order $m$ to the 
$d$-Veronese surface $V_{2,d}$. 
It is known that the tangential variety to $V_{2,d}$ is not $h$-defective 
unless $d=3$ and $h=1$ (see \cite{CGG} and \cite{B}) and that the 
$2$-osculating variety of $V_{2,d}$ is not $h$-defective unless $d=4$ and 
$h=1$ (see \cite{BF1}). Other related results are collected in the very 
recent pre-print \cite{BCGI}. Here instead we prove the following: 

\begin{Proposition}\label{special} 
Fix integers $d \ge 1$, $h \ge 1$, $m \ge 1$. 
If  
$$
\frac{(d+2)(d+1)}{2}  \le 2(h+1) + (h+1)\frac{(m+2)(m+1)}{2} 
$$ 
then $T(m,V_{2,d})$ is $h$-defective for:

(a) $h=1$, $m+2 \le d \le 2m+2$;

(b) $h=2$, $3(m+2)/2 \le d \le 2m+2$;

(c) $h=4$, $2m+4 \le d \le (5m+8)/2$;

(d) $h=5$, $12(m+2)/5 \le d \le (5m+8)/2$;

(e) $h=6$, $21(m+2)/8 \le d \le (8m+14)/3$;

(f) $h=7$, $48(m+2)/17 \le d \le (17m+32)/6$.

\noindent If instead 
$$
\frac{(d+2)(d+1)}{2}  \ge 2(h+1) + (h+1)\frac{(m+2)(m+1)}{2} 
$$ 
then $T(m,V_{2,d})$ is $h$-defective for:

(a) $h=1$, $m+1 \le d \le 2m$;

(b) $h=2$, $3(m+1)/2 \le d \le 2m$;

(c) $h=4$, $2m+2 \le d \le (5m+3)/2$;

(d) $h=5$, $12(m+1)/5 \le d \le (5m+3)/2$;

(e) $h=6$, $21(m+1)/8 \le d \le (8m+6)/3$;

(f) $h=7$, $48(m+1)/17 \le d \le (17m+15)/6$.
\end{Proposition}

\begin{Theorem}\label{nonspecial} 
Fix integers $d \ge 1$, $h \ge 1$, $m \le 18$.
Then
$$
\dim (\Sec_h(T(m,V_{2,d}))) = \expdim (\Sec_h(T(m,V_{2,d}))) =
$$
$$
= \min \left\{2(h+1) + (h+1)\frac{(m+2)(m+1)}{2} 
- 1, \frac{(d+2)(d+1)}{2} - 1\right\}
$$ 
for 

(a) $h=1$ and either $d < m+1$ or $d > 2m+2$;

(b) $h=2$ and either $d < 3(m+1)/2$ or $d > 2m+2$;

(c) $h=4$ and either $d < 2m+2$ or $d > (5m+8)/2$;

(d) $h=5$ and either $d < 12(m+1)/5$ or $d > (5m+8)/2$;

(e) $h=6$ and either $d < 21(m+1)/8$ or $d > (8m+14)/3$;

(f) $h=7$ and either $d < 48(m+1)/17$ or $d > (17m+32)/6$;

(g) $h=3$ or $h \ge 8$.

\end{Theorem}

We deduce the above statements from the classification of $(-1)$-special 
linear systems on the projective plane due to Ciliberto and Miranda 
(\cite{CilMir:00}). Indeed, our proof shows that Theorem~\ref{nonspecial} 
holds for every integer $m$ such that any special linear system of plane 
curves with base points of equal multiplicity $m+2$ is $(-1)$-special.

We are sincerely grateful to the anonymous referee for her/his crucial 
remarks on an earlier version of the present paper.

This research is part of the T.A.S.C.A. project of I.N.d.A.M., supported by 
P.A.T. (Trento) and M.I.U.R. (Italy).

\section{The results}
Let $X \subset \PP^r$ be a non-degenerate integral projective variety of 
dimension $n$ defined over the complex field $\C$. 
If $p \in X$ and $m$ is a non-negative integer, let 
$T^m_p(X)$ denote the $m$-osculating space to $X$ at $p$ (in particular, 
we have $T^0_p(X)= \{p\}$ and $T^1_p(X)=T_p(X)$, the usual tangent space).  
Fixed non-negative integers $m_0, \ldots, m_h$, define the higher order 
join
$$
J(m_0, \ldots, m_h, X) := \overline{\bigcup_{p_0, \ldots, p_h}
< T^{m_0}_{p_0}(X), \ldots, T^{m_h}_{p_h}(X) >},
$$
where $p_0, \ldots, p_h$ are general points on $X$. In particular, 
if $m_0 = \ldots = m_h = m$, we have
$$
J(m, \ldots, m, X) = \Sec_h(T(m,X)),
$$
where $T(m,X)$ is the $m$-osculating variety of $X$ (see for instance 
\cite{BF2}, Definition~1) and $\Sec_h$ denotes the $h$-secant variety 
(see for instance \cite{DioFon:01}, Definition~1.1).

\begin{Lemma}\label{Terracini} Notation as above. 

(i) The expected dimension of $J(m_0, \ldots, m_h, X)$ is 
$$
\expdim(J(m_0, \ldots, m_h, X)) = \min \left\{ 
(h+1)n + \sum_{i=0}^h {{m_i + n} \choose n} - 1, r \right\}.
$$

(ii) There is a natural inclusion 
$$
T_P(J(m_0, \ldots, m_h, X)) \subseteq 
<T^{m_0+1}_{p_0}(X), \ldots, T^{m_h+1}_{p_h}(X) >,
$$
where $p_0, \ldots, p_h$ are general points on $X$ and $P$ is general in 
$< T^{m_0}_{p_0}(X), \ldots,$ $T^{m_h}_{p_h}(X) >$.

(iii) If $r \ge (h+1)n + \sum_{i=0}^h {{m_i + n} \choose n} - 1$ and 
$\dim < T^{m_0}_{p_0}(X), \ldots, T^{m_h}_{p_h}(X) > = 
\sum_{i=0}^h {{m_i + n} \choose n} - 1 - \delta$, then 
$$
\dim J(m_0, \ldots, m_h, X)) \le (h+1)n + \sum_{i=0}^h {{m_i + n} 
\choose n} - 1 - \delta. 
$$

(iv) If $r \ge \sum_{i=0}^h {{m_i + n+1} \choose n} - 1$ and  
$\dim < T^{m_0+1}_{p_0}(X), \ldots, T^{m_h+1}_{p_h}(X) > =
\sum_{i=0}^h {{m_i + n+1} \choose n} - 1$, then 
$$
\dim(J(m_0, \ldots, m_h, X)) = \expdim(J(m_0, \ldots, m_h, X)).
$$

(v) If $n = 2$, $\dim <T^{m}_{p_0}(X), \ldots, T^{m}_{p_h}(X) > = 
\min \left\{(h+1) \frac{(m+2)(m+1)}{2} - 1, r \right\}$ and 
$\dim \Sec_h(T(m,X)) < \min \left\{2(h+1) + (h+1)\frac{(m+2)(m+1)}{2} - 1, r 
\right\}$, then 
$$
T_P(\Sec_h(T(m,X))) = <T^{m+1}_{p_0}(X), \ldots, T^{m+1}_{p_h}(X) >,
$$ 
where $p_0, \ldots, p_h$ are general points on $X$ and $P$ is general in 
$< T^{m}_{p_0}(X), \ldots,$ $T^{m}_{p_h}(X) >$.
\end{Lemma} 

\proof For $0 \le i \le h$ let 
\begin{eqnarray*}
p_i: U_i \subseteq \C^n &\longrightarrow& X \\
t_{i1}, \ldots, t_{in} &\longmapsto& p_i(t_{i1}, \ldots, t_{in})
\end{eqnarray*}
be a local parametrization of $X$ (in the euclidean topology)
centered in a general point of $X$. 
A general point of $J(m_0, \ldots, m_h, X)$ is of the form: 
$$ 
P = p_0(t_{01}, \ldots, t_{0n}) + 
\sum_{i=1}^h \lambda_i p_i(t_{i1}, \ldots, t_{in}) + 
\sum_{ {i = 0} \atop{ {1 \le \vert I 
\vert \le m_i}} }^h \lambda_i^I p_i^I(t_{i1}, \ldots, t_{in}) 
$$
where $p^I$ denotes the derivative of $p$ corresponding to the 
multi-index $I = (a_1, a_2, \ldots)$ with every $a_j \in \{1, \ldots, n \}$.
Notice that $P$ depends on $(h+1)n + \sum_{i=0}^h {{m_i + n} \choose n} - 1$ 
parameters, hence (i) follows. 

The tangent space $T_P(J(m_0, \ldots, m_h))$ is the linear span 
$$ 
< \{ P \} \cup \left\{ \frac{\partial P}{\partial t_{ij}} 
\right\}_{ {i = 0, \ldots, h} \atop {j = 1, \ldots, n}}  
\cup \left\{ \frac{\partial P}{\partial \lambda_i} 
\right\}_{i = 0, \ldots, h} 
\cup \left\{ \frac{\partial P}{\partial \lambda^I_i} 
\right\}_{ {i = 0, \ldots, h}  
\atop { 1 \le \vert I \vert \le m_i }} > =  
$$
$$
= < \{ p_i \}_{i = 0, \ldots, h} \cup \{ p_i^I \}_{ {i = 0, \ldots, h} 
\atop { 1 \le \vert I \vert \le m_i }} \cup 
\left\{ \sum_{ {\vert I \vert = m_i} } 
\lambda_i^I p_i^{I \cup \{j\}} \right\}_{ {i = 0, \ldots, h} \atop 
{j = 1, \ldots, n}} >,
$$
hence (ii), (iii), and (iv) follow. 

Finally we have to check (v). Since 
$<T^{m}_{p_0}(X), \ldots, T^{m}_{p_h}(X) >$ is of the expected dimension 
but $\Sec_h(T(m,X))$ is not, there is at least one $i \in \{0, \ldots, h \}$ 
and $j(i) \in \{1,2\}$ such that $\sum_I \lambda_i^I p_i^{I \cup \{j(i)\}}$ is 
a linear combination of the other points spanning $T_P(\Sec_h(T(m,X))$. 
Since the coefficients $\lambda_i^I$ are general, by specializing all but 
one coefficients $\lambda_i^I$ to $0$ and the remaining one to $1$ in all 
possible ways, we obtain:
\begin{eqnarray*}
T_P(\Sec_h(T(m,X)) &=&
< \{ p_k \}_{k = 0, \ldots, h} \cup \{ p_k^I \}_{ {k = 0, \ldots, h} 
\atop { 1 \le \vert I \vert \le m }} \cup \\ 
& & \cup \{ p_i^I \}_{\vert I \vert = m+1} \cup 
\left\{ \sum_{ {\vert I \vert = m} } 
\lambda_k^I p_k^{I \cup \{j\}} \right\}_{ {k \ne i} \atop 
{j = 1, \ldots, n}} >.
\end{eqnarray*}
Moreover, since $m_0 = \ldots = m_h = m$, the same is true for every 
$i \in \{0, \ldots, h \}$, hence the claim follows. 

\qed

Now we turn to the promised applications. 
We point out that our arguments rely on the main classification result 
of \cite{CilMir:00}, from which we borrow also notation and terminology.

\emph{Proof of Proposition~\ref{special}.}
In the former case, we have 
$$
\expdim (\Sec_h(T(m,V_{2,d}))) = \frac{(d+1)(d+2)}{2} - 1
$$
for every $m \ge 1$, hence the claim follows directly from 
Lemma~\ref{Terracini}~(ii) and the speciality of the  
linear system $\mathcal{L}_d((m+2)^{h+1})$
(see \cite{CCMO:03}, Theorem~2.4). 
In the latter case, the claim follows from Lemma~\ref{Terracini}~(iii) 
and the speciality of the linear system $\mathcal{L}_d((m+1)^{h+1})$
(see again \cite{CCMO:03}, Theorem~2.4). 

\qed 

\emph{Proof of Theorem~\ref{nonspecial}.}
Notice first of all that, since $m+2 \le 20$, any special 
$\mathcal{L}_d((m+2)^{h+1})$ is $(-1)$-special by \cite{CCMO:03}.
Next, if $\frac{(d+1)(d+2)}{2} \ge (h+1)\frac{(m+3)(m+2)}{2}$, 
just apply Lemma~\ref{Terracini}~(iv). 
Finally, if $\frac{(d+1)(d+2)}{2} \le (h+1)\frac{(m+3)(m+2)}{2}$, 
argue by contradiction and conclude by Lemma~\ref{Terracini}~(v).   

\qed

\noindent
Edoardo Ballico \newline
Universit\`a degli Studi di Trento \newline
Dipartimento di Matematica \newline
Via Sommarive 14 \newline
38050 Povo, Trento, Italy \newline
e-mail: ballico@science.unitn.it

\vspace{0.5cm}

\noindent
Claudio Fontanari \newline
Universit\`a degli Studi di Trento \newline
Dipartimento di Matematica \newline
Via Sommarive 14 \newline
38050 Povo, Trento, Italy \newline
e-mail: fontanar@science.unitn.it


\begin{thebibliography}{99}

\bibitem{Adl:87} B.~\r{A}dlandsvik: Joins and higher secant varieties, 
Math. Scand. 62 (1987), 213--222.

\bibitem{B} E. Ballico, On the secant varieties to the tangent developable 
of a Veronese variety. Pre-print (2003).

\bibitem{BF1} E.~Ballico and C.~Fontanari: On the secant varieties to the 
osculating variety of a Veronese surface. Cent. Eur. J. Math. 1 (2003), 
315--326. 

\bibitem{BF2} E.~Ballico and C.~Fontanari: On a Lemma of Bompiani. 
Rend. Sem. Mat. Univ. Politec. Torino (to appear).

\bibitem{BF3} E.~Ballico and C.~Fontanari: Birational geometry of defective 
varieties. Pre-print (2003).

\bibitem{BCGI} A.~Bernardi, M.~V.~Catalisano, A.~Gimigliano, M~Id\`a: 
Osculating varieties of Veronesean and their higher secant varieties. 
Pre-print math.AG/04033132 (2004). 

\bibitem{CGG} M. V. Catalisano, A. V. Geramita, A.~Gimigliano: On the 
secant variety to the tangential varieties of a
Veronesean, Proc. Amer. Math. Soc. 130 (2001), no. 4, 975--985.

\bibitem{CCMO:03}C.~Ciliberto, F.~Cioffi, R.~Miranda, F.~Orecchia:
Bivariate Hermite Interpolation and Linear Systems of Plane Curves with
Base Fat Points, Proceedings ASCM 2003, Lecture Notes Series on Computing 10,
World Scientific Publ., Singapore/River Edge, USA (2003), 87--102.

\bibitem{CilMir:00} C.~Ciliberto and R.~Miranda: Linear systems of plane 
curves with base points of equal multiplicity. Trans. Amer. Math. Soc. 352 
(2000), no. 9, 4037--4050.

\bibitem{Dale:84} M.~Dale: Terracini's lemma and the secant
variety of a curve. Proc. London Math. Soc. (3), 49 (1984), 329-339.

\bibitem{DioFon:01} C.~Dionisi and C.~Fontanari: 
Grassmann defectivity \`a la Terracini. 
Matematiche (Catania) 56 (2001), 245--255.

\bibitem{Terracini:11} A.~Terracini: Sulle $V_k$ per cui la
variet\`a degli $S_h$-$h+1$ seganti ha dimensione minore dell'ordinario.
Rend. Circ. Mat. Palermo 31 (1911), 392-396.

\bibitem{Terracini:22} A.~Terracini: Sulle superficie i cui spazi osculatori 
presentano particolari incidenze coi piani tangenti o fra loro. 
Atti Soc. Natur. e Matem. Modena, V, 6 (1922), 34-58.

\end{thebibliography}
\end{document}